\documentclass[a4paper]{article}
\pdfoutput=1
\usepackage{titlesec}

\titleformat{\section}
{\normalfont\fontsize{16}{16}\bfseries}{\thesection}{1em}{}
\titleformat{\subsection}
{\normalfont\fontsize{14}{14}\bfseries}{\thesubsection}{1em}{}
\titleformat{\subsubsection}
{\normalfont\fontsize{12}{12}\bfseries}{\thesubsubsection}{1em}{}

\usepackage[font=small,labelfont=bf]{caption}

\usepackage[toc,page]{appendix}

\usepackage[a4paper,
left=2cm,
right=2cm,
top=2cm,
bottom=2cm]{geometry}
\usepackage{setspace}
\usepackage{amsmath, nccmath}
\usepackage{amssymb}
\usepackage{mathtools}
\usepackage{graphicx}
\usepackage{indentfirst}
\usepackage[numbers]{natbib}

\numberwithin{equation}{section}
\graphicspath{ {./graphs/} }
\usepackage{fancyhdr}
\makeatletter
\renewcommand\footnotesize{%
	\@setfontsize\footnotesize\@ixpt{11}%
	\abovedisplayskip 8\p@ \@plus2\p@ \@minus4\p@
	\abovedisplayshortskip \z@ \@plus\p@
	\belowdisplayshortskip 4\p@ \@plus2\p@ \@minus2\p@
	\def\@listi{\leftmargin\leftmargini
		\topsep 4\p@ \@plus2\p@ \@minus2\p@
		\parsep 2\p@ \@plus\p@ \@minus\p@
		\itemsep \parsep}%
	\belowdisplayskip \abovedisplayskip
}
\makeatother

\begin{document}
		\title{Methods for the Numerical Analysis of Boundary Value Problem of Partial Differential Equations Based on Kolmogorov Superposition Theorem}
		\date{}
		
		\author{Korney Tomashchuk}
		\vspace{1cm}
		
		\maketitle		

	\begin{abstract}
		\noindent
		This research introduces a new method for the transition from partial to ordinary differential equations
		that is based on the Kolmogorov superposition theorem. In this paper, we discuss the numerical implementation of the Kolmogorov theorem and propose an approach that allows us to apply the theorem to represent partial derivatives of multivariate function as a combination of ordinary derivatives of univariate functions. We tested the method by running a numerical experiment with the Poisson equation. As a result, we managed to get a system of ordinary differential equations whose solution
		coincides with a solution of the initial partial differential equation.
	\end{abstract}
	
	\vspace{1cm}
	\noindent
	\textbf{Keywords:} Kolmogorov superposition theorem, partial differential equations, boundary value problem of ordinary differential equations.

	\newpage
	\onehalfspacing
	\section{Introduction}\label{Section Intro}
		Partial differential equations are a powerful instrument of mathematics that allows to model various processes of different complexities. PDEs\footnote{Partial differential equations} are applied to many real-life problems. For example, this family of equations may be used for modeling blood flow in the human venous system or for price evaluation of some financial instrument. However, very often partial differential equations are very hard to solve. Some of them we even can not solve at all, so we use modern methods of machine learning to approximate the solution. Obviously, bifurcation analysis of partial differential equations is even more complex task to solve and mathematicians have not managed to develop a general procedure for how to do it.
		
		Ordinary differential equation is a much simpler class of equations. The number of processes that can be described by them is limited, but we can solve any type of them using methods of numerical integration and analyse behaviour of a dynamic system generated by them using bifurcation analysis.
		
		Thus in this research, we worked on the derivation of the method for the transition from partial differential equations to equivalent ordinary differential equations, as this will allow us to apply all the existing methods of solution and analysis of ODEs\footnote{Ordinary differential equation} to PDEs.
		
		The main research aim of our work is to develop a method, based on the Kolmogorov superposition theorem, for the transition from partial to ordinary differential equations, such that we will be able not only to solve PDE's using method of simple numerical integration, but we will be able to apply existing methods of bifurcation analysis to partial differential equations.
	
	\section{Problem Statement}
	The problem can be formulated as follows. Consider a boundary problem of partial differential equations that is given by (\ref{PDE general}) and corresponding boundary conditions (\ref{BC general})
	
	\begin{equation}\label{PDE general}
		F\Big(x_1, \dots, x_n, f(x_1, \dots, x_n), \frac{\partial f(x_1,\dots,x_n)}{\partial x_{i_1}}, \frac{\partial^2 f(x_1,\dots,x_n)}{\partial x_{i_1} \partial x_{i_2}}, \dots, \frac{\partial^N f(x_1,\dots,x_n)}{\partial x_{i_1} \dots \partial x_{i_N}}\Big) = 0
	\end{equation}

	\begin{equation}\label{BC general}
		\begin{split}
			\varphi_{k_1}\Big(x_1, \dots, x_n, f(x_1, \dots, x_n), \frac{\partial f(x_1,\dots,x_n)}{\partial x_{i_1}},& \frac{\partial^2 f(x_1,\dots,x_n)}{\partial x_{i_1} \partial x_{i_2}}, \\
			& \dots, \frac{\partial^{N-1} f(x_1,\dots,x_n)}{\partial x_{i_1} \dots \partial x_{i_N}}\Big)\Bigg|_{x_{t_1} = x_{t_1}^{min}} = 0 \\
			\varphi_{k_2}\Big(x_1, \dots, x_n, f(x_1, \dots, x_n), \frac{\partial f(x_1,\dots,x_n)}{\partial x_{i_1}},& \frac{\partial^2 f(x_1,\dots,x_n)}{\partial x_{i_1} \partial x_{i_2}}, \\
			& \dots, \frac{\partial^{N-1} f(x_1,\dots,x_n)}{\partial x_{i_1} \dots \partial x_{i_N}}\Big)\Bigg|_{x_{t_2} = x_{t_2}^{max}} = 0 \\
		\end{split}
	\end{equation}
	where $k_1,k_2 \in \mathbb{N}, t_1,t_2\in\{1,\dots,n\}, i_{j} \in \{1,\dots,n\}$, $j = 1,\dots,N$.
	We seek to find ordinary differential equation (\ref{ODE general})
	\begin{equation}\label{ODE general}
		\tilde{F}\Big(x, y(x), \frac{dy(x)}{dx}, \frac{d^2y(x)}{dx^2}, \dots, \frac{d^{\tilde{N}}y(x)}{dx^{\tilde{N}}}\Big) = 0,
	\end{equation}
	such that solution of the equation coincides with the solution of initial problem (\ref{PDE general}, \ref{BC general}). To solve this kind of problems we developed a method based on KST\footnote{Kolmogorov superposition theorem} that we present in the sections 4-7.
	
	\section{Related Works}
	\subsection{Method for transition from partial to ordinary differential equations}
	As partial differential equations have ben used for mathematical modeling in various fields for a long time, the idea of reducing them to ordinary differential equations have been studied by many mathematicians. In a review paper \cite{p2oReview} group of researchers give a detailed review of existing methods for reduction of PDEs to ODEs. They carry various experiments in order to compare and outline pros and cons of each method. Two other mathematicians in their work \cite{singhatanadgid2019kantorovich} write about superiority of the Kontorovich method on other existing methods supporting their opinions with several numerical experiments. However, existing methods have their own disadvantages. Some of them use approximations that lead to uncertainties. Other methods that use orthogonal basis, such as Fourier based methods, break bifurcation picture of the equations making conclusions of bifurcation theory irrelevant. Our method is based on a strict equality that is known as Kolmorogov superposition theorem.
	\subsection{Kolmogorov superposition theorem}
	The method that we propose in this paper relies heavily on Sprecher's formulation of the Kolmogorov superposition theory that he introduced in the article \cite{sprecher1996numerical}. This formulation is more convenient in terms of numerical realisation. In the same article he proposes algorithms for calculating inner and outer functions that are required in order to apply the theorem on practice.
	
	Despite that Sprecher's formulation of the theorem was correct, the algorithm for calculating internal functions was wrong. For certain values of parameters the monotonicity property of these functions is violated. Another Mathematician Mario Koppen who also dedicated part of his work to study the Kolmogorov theorem, spotted this mistake and corrected. Proper algorithm for constructing inner functions was introduced in the paper \cite{koppen2002training}.
	
	Jonas Actor in his PhD thesis \cite{actor2018computation} on the study of the properties of the Kolmogorov theorem, give a review of existing formulations of Kolmogorov theorem, and proposes a new type of improved inner functions that are Lipschitz continuous. Initial functions that were proposed by David Sprecher have very steep slope on some segments. This can lead to computational problems. In contrast, Lipschitz continuous functions have more controlled slope.
	
	T. Hedberb also studied properties of the Komogorov superposition theorem. In the work \cite{hedberg1971kolmogorov} he derived his own formulation of the theorem that does not use strictly defined inner functions. Instead, he have showed his reformulation is correct for quasi-all functions if only they satisfy certain conditions described in \cite{hedberg1971kolmogorov}.
	
	\section{Numerical implementation of Kolmogorov superposition theorem}
	\subsection{Theorem formulation review}
	The method for transition from partial to ordinary differential equations is based on Sprecher's reformulation of Kolmogorov theorem. We start with providing a quick overview of the theorem. Consider a function of $n$ variables $f(x_1,\dots,x_n)$. As stated in the paper \cite{sprecher1996numerical}, $f(x_1,\dots,x_n)$ can be written as a sum of functions of one variable (\ref{KST})
	
	\begin{equation}\label{KST}
		f(x_1, x_2, \dots, x_n) = \sum_{q=0}^{2n}\Phi_q\Big[\sum_{p=1}^{n}\alpha_p\psi(x_p + aq)\Big],
	\end{equation}
	where constants $a, \alpha_p, p = 1,\dots,n$ and inner function $\psi(x)$ are the same for different $f(x_1,\dots,x_n)$. Moreover, function $\psi(x)$ is a specific function that can be calculated numerically using algorithm of Mario Koppen described in \cite{koppen2002training}. Initial algorithm proposed by David A. Sprecher contained a mistake for some combination of algorithm parameters. Further, we provide a quick review of a correct algorithm.
	
	\subsection{Numerical implementation of inner $\psi(x)$ function}
	Let's introduce notation $D$ for the set of terminating rational numbers that are defined according to equation (\ref{dk}).
	
	\begin{equation}\label{dk}
		d_k = \sum_{j=1}^{k}\frac{i_j}{\gamma^k}, i_j = 0, \dots, \gamma-1,
	\end{equation}	
	where $k, \gamma$ are parameters of the algorithm such that $\gamma \ge 2n+2$ and $k \in \mathbb{N}$. Next, we define constants $a, \alpha_p$ according to the formulas (\ref{a alpha}).
	
	\begin{equation}\label{a alpha}
		\begin{split}
			a &= \frac{1}{\gamma(\gamma - 1)} \\
			\alpha_p &= \sum_{j=1}^{\infty}\gamma^{-(p-1)\frac{n^{p}-1}{n-1}}
		\end{split}
	\end{equation}
	Then, the inner function $\psi(x)$ is recursively defined by the formula (\ref{psi})
	\begin{equation}\label{psi}
		\psi^k(d_k) = \begin{cases}
			d_k, \text{ for } k=1\\
			\psi^{k-1}(d_k - \frac{i_k}{\gamma_k}) + \frac{i_k}{\gamma^{(n^k-1)/(n-1)}} \text{ for } k>1, i_k<\gamma-1\\
			\frac{1}{2} \Big(\psi_k(d_k - \frac{1}{\gamma^k}) + \psi_{k-1}(d_k + \frac{1}{\gamma_k})\Big) \text{ for } k>1, i_k=\gamma-1
		\end{cases}
	\end{equation}
	where $i_k$ are decimal integers of rational number $d_k = 0.i_1i_2\dots i_k$. In most of the problems the function $f(x_1,\dots,x_n)$ is defined beyond the unit cube. For $x > 1$, function $\psi(x)$ is calculated using formula (\ref{psi > 1}).
	\begin{equation}\label{psi > 1}
		\psi(x) = \psi(x - [x]) + [x],
	\end{equation}
	where $[x]$ is the integer part of $x$.
	
	According to the work \cite{koppen2002training} of Mario Koppen  $\gamma = 10$ is the best value of parameter for $n=2$. Graphs for function $\psi(x)$ for $k=2,3,4$ and $\gamma=10$ are depicted in the \textbf{Figure 1}. Function $\psi(x)$ is defined only at certain points $d_k$ but for the better representation we use linear interpolation.
	\begin{figure}[h]\label{Fig1 psi}
		{\includegraphics[width=6in, height=2in]{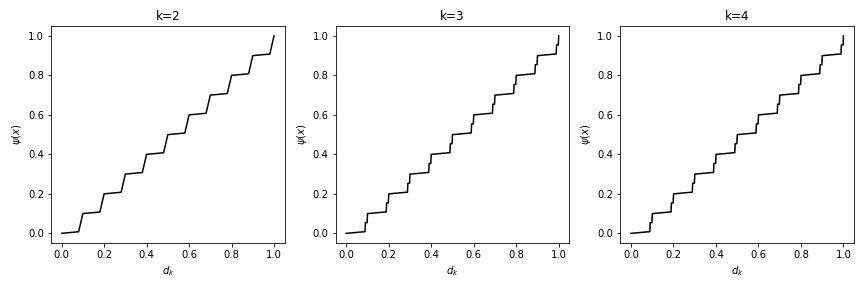}}
		\centering\caption{$\psi(x)$ for $k=2,3,4$ and $\gamma=10$}
	\end{figure}

\begin{figure}[h]
		{\includegraphics[width=4in, height=1.5in]{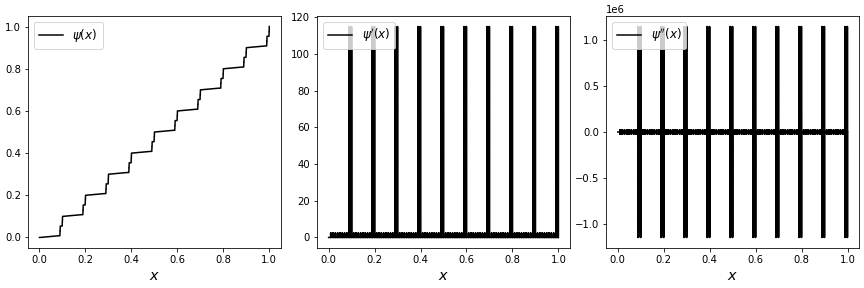}}
		\centering\caption{Graphs of $\psi(x), \psi'(x), \psi''(x)$ with parameters $k=4, \gamma=10$.}
		\label{Fig2 psi ders}
	\end{figure}

	As we are going to apply Kolmogorov theorem to the PDE problems, we should be able to define derivatives of inner functions. $\psi(x)$ cannot be analytically differentiated. Thus, we consider difference analogue of the derivative (\ref{psi der}).
	
	\begin{equation}\label{psi der}
		\psi'(x) = \frac{\psi(x + \Delta) - \psi(x)}{\Delta},
	\end{equation}
	where $\Delta$ is an increment that depends on the choice of the parameters $\gamma$ and $k$. As $\psi(x)$ is defined only at a finite number of points $d_k$, we take $\Delta = \gamma^{-k}$ equal to the size of the step between two neighbour points from the set $D_k$. Graphs of the first two derivates of the function $\psi(x)$ with parameters $k=4, \gamma=10$ are depicted in the \textbf{Figure \ref{Fig2 psi ders}}.
	
	\section{Representation of a multivariate function as a sum of functions of one variable}
	Transition from partial to ordinary differential equations consists of two main parts. First, we should represent function of several variable as a combination of functions of one variable. Second, substitute partial derivatives with corresponding ordinary. To complete the first stage, we use Kolmogorov theorem and the fact that $a = \frac{1}{\gamma(\gamma-1)}$ is a small parameter. For recommended value of $\gamma = 10$ for $n=2$, $a$ is equal to $\frac{1}{90}$. As there are no upper bound for a value of $\gamma$, we can consider taking $\gamma \rightarrow \infty$, as a result $a\rightarrow 0$. Thus, we can consider taking Taylor series of (\ref{KST}) in the vicinity $a = 0$ (see equation (\ref{Taylor KST})).
	
	\begin{equation}\label{Taylor KST}
		\begin{split}
			\sum_{q=0}^{2n}\Phi_q\Big(\sum_{p=0}^{n}&\alpha_p\psi(x_p + qa)\Big) = \sum_{m=0}^{M}\frac{a^m}{m!}\frac{d^m}{da^m}\Bigg(\sum_{q=0}^{2n}\Phi_q\Big(\sum_{i=p}^{n}\alpha_p\psi(x_p + aq)\Big)\Bigg)\Bigg|_{a=0}
		\end{split}
	\end{equation}
	
	In order to get representation of Taylor series with $M$ terms of KST theorem, we use Bell's polynomials and Faà di Bruno's formula. Definition of both are provided below.

	\subsection{Bell's polynomial}
	
	By definition, Bell's polynomials are defined by the formula 
	
	\begin{equation}\label{Bell Poly}
		B_{n,k}(x_1, \dots, x_{n-k+1}) = \sum\bigg(\frac{n!}{j_1!j_2!\dots j_{n-k+1}!}\prod_{i=1}^{n-k+1}\Big(\frac{x_i}{i!}\Big)^{j_i}\bigg),
	\end{equation}
	where outer sum is taken over all non-negative numbers $j_1,j_2,\dots,j_{n-k+1}$ that satisfy two conditions (\ref{Bell j Cond 1}, \ref{Bell j Cond 2}).
	\begin{equation}\label{Bell j Cond 1}
		\sum_{i=1}^{n-k+1}j_{i} = k
	\end{equation}
	\begin{equation}\label{Bell j Cond 2}
		\sum_{i=1}^{n-k+1}ij_{i} = n
	\end{equation}
	
	\subsection{Faà di Bruno's formula}
	
	Faà di Bruno's formula is given by the equation (\ref{Faa-di-Bruno})
	\begin{equation}\label{Faa-di-Bruno}
		\frac{d^m}{dx^m}f(g(x)) = \sum\bigg(\frac{m!}{j_1!j_2!\dots j_m!}\cdot f^{(j_1+\dots+j_m)}\big(g(x)\big)\prod_{i=1}^{m}\Big(\frac{g^{(i)}(x)}{i!}\Big)^{j_i}\bigg),
	\end{equation}
	where external sum is taken over all nonnegative integers that satisfy second condition of Bell's polynomial (\ref{Bell j Cond 2}). It is used to define high-order derivatives of complex functions. Using Bell's polynomials, we can rewrite formula (\ref{Faa-di-Bruno}) in the following form.
	
	\begin{equation}\label{Faa-di-Bruno Bell}
		\frac{d^m}{dx^m}f(g(x)) = \sum_{k=0}^{m}f^{(k)}\big(g(x)\big)B_{m,k}\big(g'(x), g''(x),\dots,g^{(m-k-1)}(x)\big)
	\end{equation}
	
	\subsection{Taylor series of Kolmogorov superposition theorem}
	
	We can substitute formula (\ref{Faa-di-Bruno Bell}) to achieve compact representation of Taylor series (\ref{Taylor KST}) as shown in the equation (\ref{der Bell})
	
	\begin{equation}\label{der Bell}
		\begin{split}
			&\frac{d^m}{da^m}\Bigg(\sum_{q=0}^{2n}\Phi_q\Big(\sum_{p=1}^{n}\alpha_p\psi(x_p + aq)\Big)\Bigg)\Bigg|_{a=0} =	\sum_{k=0}^{m}\sum_{q=0}^{2n}\Phi_q^{(k)}\Big(\sum_{p=1}^{n}\alpha_p\psi(x_p)\Big)\cdot\\
			&\cdot B_{m,k}\big(q\sum_{p=1}^{n}\alpha_p\psi'(x_p), q^2\sum_{p=1}^{n}\alpha_p\psi''(x_p), \dots, q^{(m-k-1)}\sum_{p=1}^{n}\alpha_p\psi^{(m-k-1)}(x_p)\big)
		\end{split}
	\end{equation}
	Using definition of Bell's polynomial and applying sequence of simple transformations, we get more convinient representation (\ref{KST Bell short}).
	
	\begin{equation}\label{KST Bell short}
		\begin{split}
			&B_{m,k}\big(q\sum_{p=1}^{n}\alpha_p\psi'(x_p), q^2\sum_{p=1}^{n}\alpha_p\psi''(x_p), \dots, q^{(m-k-1)}\sum_{p=1}^{n}\alpha_p\psi^{(m-k-1)}(x_p)\big) =\\
			&= \sum\bigg(\frac{m!}{j_1!j_2!\dots j_{m-k+1}!}\prod_{i=1}^{m-k+1}\Big(\frac{q^i\sum_{p=1}^{n}\alpha_p\psi^{(i)}(x_p)}{i!}\Big)^{j_i}\bigg) = q^{m}\tilde{B}_{m,k}(x_1,x_2\dots,x_n)
		\end{split}
	\end{equation}
	where
	\begin{equation}\label{Bell tilde}
		\tilde{B}_{m,k}(x_1,x_2\dots,x_n) = \sum\bigg(\frac{m!}{j_1!j_2!\dots j_{m-k+1}!}\prod_{i=1}^{m-k+1}\Big(\frac{\sum_{p=1}^{n}\alpha_p\psi^{(i)}(x_p)}{i!}\Big)^{j_i}\bigg)
	\end{equation}
	Then, Taylor series of KST representation of function of $n$ variables with arbitrary number of terms ($M$) can be represented with the equation (\ref{Taylor KST Bell raw})
	\begin{equation}\label{Taylor KST Bell raw}
		\begin{split}
			\sum_{q=0}^{2n}\Phi_q\Big(\sum_{p=0}^{n}\alpha_p\psi(x_p + qa)\Big) = \sum_{m=0}^{M}\sum_{k=0}^{m}\tilde{B}_{m,k}(x_1,\dots,x_n)\sum_{q=0}^{2n}\frac{a^mq^m}{m!}\Phi_q^{(k)}\Big(\sum_{p=1}^{n}\alpha_p\psi(x_p)\Big)
		\end{split}
	\end{equation}
	Finally, by making a substitution $z = \sum_{p=1}^{n}\alpha_p\psi(x_p)$, we derive representation of multivariate $f(x_1, x_2, \dots, x_n)$ as a sum of univariate functions $\Phi_q(z)$ (see equation (\ref{Taylor KST Bell})).
	\begin{equation}\label{Taylor KST Bell}
		f(x_1, x_2, \dots, x_n) = \sum_{m=0}^{M}\sum_{k=0}^{m}\tilde{B}_{m,k}(x_1,x_2\dots,x_n)\sum_{q=0}^{2n}\frac{a^mq^m}{m!}\Phi_q^{(k)}(z),
	\end{equation}

	\section{Method for transition from partial to ordinary differential equations}\label{section p2o}
	Method for transition from partial to ordinary differential equations that we derived can be applied only to variational formulation of PDE problem. Let's consider functional $\Xi$
	\begin{equation}\label{Functional general}
		\Xi = \int_{x}L\Big(x_1, \dots, x_n, f(x_1, \dots, x_n), \frac{\partial f(x_1,\dots,x_n)}{\partial x_{i_1}}, \frac{\partial^2 f(x_1,\dots,x_n)}{\partial x_{i_1} \partial x_{i_2}}, \dots, \frac{\partial^N f(x_1,\dots,x_n)}{\partial x_{i_1} \dots \partial x_{i_N}}\Big)dx
	\end{equation}
	such that solution of variational problem (\ref{Functional general extr}) coincide with the solution of initial PDE problem (\ref{PDE general}, \ref{BC general}).
	\begin{equation}\label{Functional general extr}
		\Xi \longrightarrow \underset{x}{extr}
	\end{equation}
	
	In order to get rid of partial derivatives in functional, we use obtained formula (\ref{Taylor KST Bell}). First, we change variable of integrations. We parameterise $x_2, \dots, x_{n}$ and $x_1$ express through other parametrised variables as shown in (\ref{variable change general}).
	
	\begin{equation}\label{variable change general}
		\begin{split}
			x_2 &= \tilde{x}_2\\
			\dots& \dots \dots \\
			x_n &= \tilde{x}_n \\
			x_1 &= \psi^{-1}\Big(\frac{z - \sum_{p=2}^{n}\alpha_p\psi(x_p)}{\alpha_1}\Big)
		\end{split}
	\end{equation}
	Boundaries of integration change according to formulas (\ref{boundaries of integaration general})
	\begin{equation}\label{boundaries of integaration general}
		\begin{split}
			&x_1^{min} \mapsto z^{min} = \alpha_1\psi(x_1^{min}) + \alpha_2\psi(\tilde{x}_2) + \dots + \alpha_n\psi(\tilde{x}_n) \\
			&x_1^{max} \mapsto z^{max} = \alpha_1\psi(x_1^{max}) + \alpha_2\psi(\tilde{x}_2) + \dots + \alpha_n\psi(\tilde{x}_n) \\
		\end{split}
	\end{equation}
	and integration factor according to determinant of Jacobian matrix (\ref{Jacobian General}). Determinant of a diagonal matrix equal to product of diagonal elements.
	\begin{equation}\label{Jacobian General}
		\begin{split}
		&J =
		\begin{bmatrix}
			z'_{x_1} & z'_{x_2} & \dots & z'_{x_n}\\
			{x'_2}_{x_1} & {x'_2}_{x_2} &\dots & {x'_2}_{x_n} \\
			\dots & \dots &\dots & \dots \\
			{x'_n}_{x_1} & {x'_n}_{x_2} &\dots & {x'_n}_{x_n}
		\end{bmatrix} = 
		\begin{bmatrix}
			\alpha_1\psi'(x_1) & \alpha_2\psi'(x_2) & \dots & \alpha_n\psi'(x_n)\\
			0 & 1 &\dots & 0 \\
			\dots & \dots &\dots & \dots \\
			0 & 0 &\dots & 1
		\end{bmatrix}\\
		&\implies dx_1 = \frac{dz}{\alpha_1\psi'\circ\psi^{-1}\Big(\frac{z - \sum_{p=2}^{n}\alpha_p\psi(\tilde{x}_p)}{\alpha_1}\Big)}
		\end{split}	
	\end{equation}
	
	Secondly, we redefine partial derivatives of $f(x_1,\dots,x_n)$ as shown in (\ref{KST der f general})
	\begin{equation}\label{KST der f general}
		\begin{split}
			&\frac{\partial^N f(x_1,\dots,x_n)}{\partial x_{i_1} \dots \partial x_{i_N}} = \sum_{m=0}^{M}\sum_{k=0}^{m}\Bigg[\sum_{j=0}^{N}\bigg[\frac{\partial^N}{\partial x_{i_1} \dots \partial x_{i_N}}\Big(\tilde{B}_{m,k}(x)\sum_{q=0}^{2n}\frac{a^mq^m}{m!}\Phi_{q}^{(k)}(z)\Big)\bigg]\Bigg]
		\end{split}
	\end{equation}
	where $\tilde{B}_{m,k}$ and $z$ are known from KST theorem. Thus, in terms of unknown functions $\Phi_{q}(z)$ we get expression that contains only ordinary derivatives. Applying described procedure to the functional (\ref{Functional general}), we derive equivalent functional (\ref{Functional general KST}) that does not contain partial derivatives.
	\begin{equation}\label{Functional general KST}
		\begin{split}
			\tilde{\Xi} = \int_{\tilde{x}}\int_{z^{min}}^{z^{max}} L\Big(z, \tilde{x}_2, \dots, \tilde{x}_n, \Phi_0(z),\dots,\Phi_{2n}(z),&\\ \frac{d\Phi_0(z)}{dz}, \dots, \frac{d\Phi_{2n}(z)}{dz}&, \dots, \frac{d^N\Phi_0(z)}{dz^N}, \dots, \frac{d^N\Phi_{2n}(z)}{dz^N}\Big) dzd\tilde{x}
		\end{split}
	\end{equation}
	Then, initial PDE problem is reduced to a much simpler variational problem with ordinary derivatives (\ref{Functional general KST extr}). This problem can be solved using standard methods of calculus of variations.
	\begin{equation}\label{Functional general KST extr}
		\begin{split}
			\forall\{\tilde{x}_2,\dots,\tilde{x}_n\}, \: \tilde{x}_2 \in [\tilde{x}_2^{min},& \tilde{x}_2^{max}], \dots, \tilde{x}_n \in [\tilde{x}_n^{min}, \tilde{x}_n^{max}]:  \\
			&\tilde\Xi \rightarrow \underset{z}{extr}
		\end{split}
	\end{equation}

	\section{Numerical experiment with Poisson equation}
	\subsection{Equivalent variational formulation}
	In order to check our method, we decided to apply it to the partial differential equation with known solution which is the Poisson equation. Let's consider the following formulation (\ref{PDE Poisson}) with boundary conditions (\ref{BC Poisson}).
	\begin{equation}\label{PDE Poisson}
		\nabla^2u(x_1, x_2) = \sin(\pi x_1)\sin(\pi x_2)
	\end{equation}
	\begin{equation}\label{BC Poisson}
		u(0, x_2) = 0, u(1, x_2) = 0, u(x_1, 0) = 0, u(x_1, 1) = 0
	\end{equation}
	Solution for problem (\ref{PDE Poisson}, \ref{BC Poisson}) is given by equation (\ref{Sol Poisson}).
	\begin{equation}\label{Sol Poisson}
		u^*(x_1, x_2) = \frac{\sin(\pi x_1)\sin(\pi x_2)}{-2\pi^2}
	\end{equation}
	One can easily check that (\ref{Sol Poisson}) satisfies both equation and boundary conditions. As shown in Section \ref{section p2o}, we should first find equivalent variational formulation. For PDE problem (\ref{PDE Poisson}, \ref{BC Poisson}) corresponding problem of calculus of variations is given by (\ref{Functional Poisson 1}). The whole derivation procedure is shown in the appendix A.
	\begin{equation}\label{Functional Poisson 1}
		\begin{split}
			\int_{0}^{1}\int_{0}^{1} \bigg[&-\Big(\frac{\partial{u(x_1,x_2)}}{\partial{x_1}}\Big)^2 - \Big(\frac{\partial{u(x_1,x_2)}}{\partial{x_2}}\Big)^2 - 2\sin(\pi x_1)\sin(\pi x_2)u(x_1, x_2) -\\
			&- 2\frac{\partial^2 u(x_1,x_2)}{\partial x_1^2}u(x_1,x_2) - 2\frac{\partial^2 u(x_1,x_2)}{\partial x_2^2}u(x_1,x_2)\bigg]dx_1dx_2 \longrightarrow extr
		\end{split}
	\end{equation} 

	\subsection{ODE boundary problem of Poisson equation}
	In order to derive corresponding system of ordinary differential equations, we apply the procedure that we described in section \ref{section p2o}. Let's consider the simplest case when we remove all terms of Taylor series except the first one. Then, according to (\ref{Taylor KST Bell}) for $n=2, M=0$ function $u(x_1, x_2)$ can be written as (\ref{KST 0 f Poisson})
	\begin{equation}\label{KST 0 f Poisson}
		u(x_1,x_2) = \sum_{q=0}^{4}\Phi_q(z),
	\end{equation}
	where $z = \alpha_1\psi(x_1) + \alpha_2\psi(x_2)$. Further, following algorithm described above we apply substitutions (\ref{KST 0 variable substitution Poisson})
	\begin{equation}\label{KST 0 variable substitution Poisson}
		\begin{split}
			&x_2 = \tilde{x}_2\\
			&x_1 = \psi^{-1}\Big(\frac{z - \alpha_2\psi(\tilde{x}_2)}{\alpha_1}\Big)
		\end{split}
	\end{equation}
	and change factor and boundaries of integration (\ref{KST 0 boundaries}, \ref{KST 0 J}).
	\begin{equation}\label{KST 0 boundaries}
		\begin{split}
			x_1 &= 0 \rightarrow z^{min} = \alpha_2\psi(\tilde{x}_2)\\
			x_1 &= 1 \rightarrow z^{max} = \alpha_1 + \alpha_2\psi(\tilde{x}_2) \\
		\end{split}
	\end{equation}
	\begin{equation}\label{KST 0 J}
		\begin{split}
			J &= \begin{bmatrix}
				z'_{x_1} & z'_{x_2} \\ {x_2}'_{x_1} & {x_2}'_{x_2}
			\end{bmatrix} = \begin{bmatrix}
				\alpha_1\psi'(x_1) & \alpha_2\psi'(x_2) \\ 0 & 1
			\end{bmatrix} \\
			&\implies dx_1 = \frac{dz}{\alpha_1\psi'(x_1)} = \frac{dz}{\alpha_1\psi'\psi^{-1}\big(\frac{z - \alpha_2\psi(\tilde{x}_2)}{\alpha_1}\big)}
		\end{split}
	\end{equation}
	Next step is to define partial derivatives up to second through ordinary derivatives of functions $\Phi_q(z)$ (see equation (\ref{KST 0 df Poisson})).
	\begin{equation}\label{KST 0 df Poisson}
		\begin{split}
			&\frac{\partial u(x_1,x_2)}{\partial x_i} = \alpha_i\psi'(x_i)\sum_{q=0}^{4}\Phi_q'(z) \\
			&\frac{\partial^2 u(x_1,x_2)}{\partial x_i^2} = \alpha_i\psi''(x_i)\sum_{q=0}^{4}\Phi_q'(z) +  \alpha_i^2\psi'^2(x_i)\sum_{q=0}^{4}\Phi_q''(z)\\
		\end{split}
	\end{equation}
	As a result, using all above substitutions, we reformulate PDE Poisson problem as a problem of finding extremum of a functional that contains only ordinary derivatives (\ref{Functional Poisson KST 0 raw})
	
	\begin{equation}\label{Functional Poisson KST 0 raw}
		\begin{split}
			\int_{0}^{1}&\int_{z_{min}}^{z_{max}} \frac{1}{\alpha_1\psi'\psi^{-1}\big(\frac{z - \alpha_2\psi(\tilde{x}_2)}{\alpha_1}\big)}\bigg[-\Big( \alpha_1\psi'\psi^{-1}\big(\frac{z - \alpha_2\psi(\tilde{x}_2)}{\alpha_1}\big)\sum_{q=0}^{4}\Phi_q'(z) \Big)^2 -\\
			& -\Big( \alpha_2\psi'(\tilde{x}_2)\sum_{q=0}^{4}\Phi_q'(z) \Big)^2 + 2\sum_{q=0}^{4}\Phi_q(z)\cdot\Big[f\Big(\psi^{-1}\big(\frac{z - \alpha_2\psi(\tilde{x}_2)}{\alpha_1}\big), \tilde{x}_2\Big) - \\
			&- \Big( \alpha_1\psi''\psi^{-1}\big(\frac{z - \alpha_2\psi(\tilde{x}_2)}{\alpha_1}\big)\sum_{q=0}^{4}\Phi_q'(z) + \alpha_1^2\psi'^2\psi^{-1}\big(\frac{z - \alpha_2\psi(\tilde{x}_2)}{\alpha_1}\big)\sum_{q=0}^{4}\Phi_q''(z) \Big) -\\
			&- 2\Big( \alpha_2\psi''(\tilde{x}_2)\sum_{q=0}^{4}\Phi_q'(z) + \alpha_2^2\psi'^2(\tilde{x}_2)\sum_{q=0}^{4}\Phi_q''(z) \Big)\bigg]dzdx_2 \longrightarrow extr
		\end{split}
	\end{equation}
	
	In order to solve problem (\ref{Functional Poisson KST 0 raw}) we should take variational derivatives with respect to each independent variable and equate to zero. First, we should identify the number of independent variables by variating  functional with respect to each of the functions $\Phi_q(z)$ and identify rank of the matrix of coefficients of the highest derivatives. To make formulas more readable, we use notation $\psi^{-1}\big(\frac{z - \alpha_2\psi(\tilde{x}_2)}{\alpha_1}\big) = \tilde{x}_1$. Then, variational derivative with respect to $\Phi_i(z)$ is equal to (\ref{Functional Poisson KST 0 eq and bc raw}).
	
	\begin{equation}\label{Functional Poisson KST 0 eq and bc raw}
		\begin{split}
			\int_{0}^{1}&\int_{z_{min}}^{z_{max}} \bigg[\frac{f(\tilde{x}_1, \tilde{x}_2)}{\alpha_1\psi'(\tilde{x}_1)}-\frac{\alpha_1^2\psi'^2(\tilde{x}_1) + \alpha_2^2\psi'^2(\tilde{x}_2)}{\alpha_1\psi'(\tilde{x}_1)}\sum_{q=0}^{4}\Phi_q''(z)- \\
			& -\frac{\alpha_1^2\psi'^2(\tilde{x}_1)\psi''(\tilde{x}_1) - \alpha_2^2\psi'^2(\tilde{x}_2)\psi''(\tilde{x}_1)}{\alpha_1^2\psi'^3(\tilde{x}_1)}\sum_{q=0}^{4}\Phi_q'(z) -\\
			&-\Big(\frac{\alpha_1\alpha_2\psi'^2(\tilde{x}_1)\psi''(\tilde{x}_1)\psi''(\tilde{x}_2)}{\alpha_1^3\psi'^5(\tilde{x}_1)}\\\
			&-\frac{\alpha_2^2\psi'^2(\tilde{x}_2)\big(3\psi''(\tilde{x}_1)-\psi(\tilde{x}_1)\psi'''(\tilde{x}_1)\big)}{\alpha_1^3\psi'^5(\tilde{x}_1)}\Big)\sum_{q=0}^{4}\Phi_q(z) \bigg]\delta{\Phi_i(z)}dzd\tilde{x}_2-\\
			-&\int_{0}^{1}\Bigg[\bigg[\frac{\alpha_2\Big(\alpha_2\psi'^2(\tilde{x}_2)\psi''(\tilde{x}_1) + \alpha_1\psi'^2(\tilde{x}_1)\psi''(\tilde{x}_2)\Big)}{\alpha_1^2\psi'^3(\tilde{x}_1)}\bigg]\sum_{q=0}^{4}\Phi_q(z)\delta{\Phi_i(z)}\Bigg|_{z_{min}}^{z_{max}} +\\
			&+\bigg[\alpha_1\psi'(\tilde{x}_1) + \frac{\alpha_2^2\psi'^2(\tilde{x}_2)}{\alpha_1\psi'(\tilde{x}_1)}\bigg]\sum_{q=0}^{4}\Phi_q(z)\delta{\Phi_i'(z)}\Bigg|_{z_{min}}^{z_{max}}\Bigg]d\tilde{x}_2 = 0
		\end{split}
	\end{equation}
	We can see that coefficient of $\Phi_i''(z)$ does not depend on the variable of variation. This implies that rank of the matrix of coefficients is 1 and that there is only one independent variable. We define this variable according to formula (\ref{one variable})
	\begin{equation}\label{one variable}
		U(z) = \sum_{q=0}^{4}\Phi_q(z)
	\end{equation}
	Making this substitution and variating functional with respect to $U(z)$, we get (\ref{Functional Poisson KST 0 eq and bc})
	\begin{equation}\label{Functional Poisson KST 0 eq and bc}
		\begin{split}
			\int_{0}^{1}&\int_{z_{min}}^{z_{max}} \bigg[\frac{f(\tilde{x}_1, \tilde{x}_2)}{\alpha_1\psi'(\tilde{x}_1)}-\frac{\alpha_1^2\psi'^2(\tilde{x}_1) + \alpha_2^2\psi'^2(\tilde{x}_2)}{\alpha_1\psi'(\tilde{x}_1)}U''(z)- \\
			& -\frac{\alpha_1^2\psi'^2(\tilde{x}_1)\psi''(\tilde{x}_1) - \alpha_2^2\psi'^2(\tilde{x}_2)\psi''(\tilde{x}_1)}{\alpha_1^2\psi'^3(\tilde{x}_1)}U'(z)\\
			&-\frac{\alpha_1\alpha_2\psi'^2(\tilde{x}_1)\psi''(\tilde{x}_1)\psi''(\tilde{x}_2) + \alpha_2^2\psi'^2(\tilde{x}_2)\big(3\psi''(\tilde{x}_1)-\psi(\tilde{x}_1)\psi'''(\tilde{x}_1)\big)}{\alpha_1^3\psi'^5(\tilde{x}_1)}U(z) \bigg]\delta{\Phi_i(z)}dzd\tilde{x}_2-\\
			-&\int_{0}^{1}\Bigg[\bigg[\frac{\alpha_2\Big(\alpha_2\psi'^2(\tilde{x}_2)\psi''(\tilde{x}_1) + \alpha_1\psi'^2(\tilde{x}_1)\psi''(\tilde{x}_2)\Big)}{\alpha_1^2\psi'^3(\tilde{x}_1)}\bigg]U(z)\delta{U(z)}\Bigg|_{z_{min}}^{z_{max}} +\\
			&+\bigg[\alpha_1\psi'(\tilde{x}_1) + \frac{\alpha_2^2\psi'^2(\tilde{x}_2)}{\alpha_1\psi'(\tilde{x}_1)}\bigg]U(z)\delta{U'(z)}\Bigg|_{z_{min}}^{z_{max}}\Bigg]d\tilde{x}_2 = 0
		\end{split}
	\end{equation}
	Variation of functional is equal to zero if and only if integrand and boundary conditions are equal to zero. Then, we get ordinary differential equation of second degree (\ref{DE2 Poisson KST 0})
	\begin{equation}\label{DE2 Poisson KST 0}
		\begin{split}
			&\frac{\alpha_1^2\psi'^2(\tilde{x}_1) + \alpha_2^2\psi'^2(\tilde{x}_2)}{\alpha_1\psi'(\tilde{x}_1)}U''(z) +\frac{\alpha_1^2\psi'^2(\tilde{x}_1)\psi''(\tilde{x}_1) - \alpha_2^2\psi'^2(\tilde{x}_2)\psi''(\tilde{x}_1)}{\alpha_1^2\psi'^3(\tilde{x}_1)}U'(z)+\\
			&+\frac{\alpha_1\alpha_2\psi'^2(\tilde{x}_1)\psi''(\tilde{x}_1)\psi''(\tilde{x}_2) + \alpha_2^2\psi'^2(\tilde{x}_2)\big(3\psi''(\tilde{x}_1)-\psi(\tilde{x}_1)\psi'''(\tilde{x}_1)\big)}{\alpha_1^3\psi'^5(\tilde{x}_1)}U(z) -\\
			&-\frac{f(\tilde{x}_1, \tilde{x}_2)}{\alpha_1\psi'(\tilde{x}_1)} = 0
		\end{split}
	\end{equation}
	with boundary conditions (\ref{BC Poisson KST 0}).
	\begin{equation}\label{BC Poisson KST 0}
		\begin{split}
			&\bigg[\frac{\alpha_2^2\psi'^2(\tilde{x}_2)\psi''(0) + \alpha_1\alpha_2\psi'^2(0)\psi''(\tilde{x}_2)}{\alpha_1^2\psi'^3(0)} + \alpha_1\psi'(0) + \frac{\alpha_2^2\psi'^2(\tilde{x}_2)}{\alpha_1\psi'(0)}\bigg]U(\alpha_2\psi(\tilde{x}_2)) = 0 \\
			&\\
			&\bigg[\frac{\alpha_2^2\psi'^2(\tilde{x}_2)\psi''(1) + \alpha_1\alpha_2\psi'^2(1)\psi''(\tilde{x}_2)}{\alpha_1^2\psi'^3(1)} + \alpha_1\psi'(1) + \frac{\alpha_2^2\psi'^2(\tilde{x}_2)}{\alpha_1\psi'(1)}\bigg]U(\alpha_1+\alpha_2\psi(\tilde{x}_2)) = 0 \\
		\end{split}
	\end{equation}
	Making substitution $U'(z) = W(z)$, we get boundary problem of ordinary differential equations of the first degree (\ref{System Poisson KST 0}) that can be solved using standard methods of numerical integration.
	\begin{equation}\label{System Poisson KST 0}
		\begin{cases}
			U'(z) = W(z) \\\\
			\begin{split}
				W'(z) &= \\
				&=\frac{f(\tilde{x}_1, \tilde{x}_2)}{\alpha_1^2\psi'^2(\tilde{x}_1) + \alpha_2^2\psi'^2(\tilde{x}_2)} -\frac{\alpha_1^2\psi'^2(\tilde{x}_1)\psi''(\tilde{x}_1) - \alpha_2^2\psi'^2(\tilde{x}_2)\psi''(\tilde{x}_1)}{\alpha_1\psi'^2(\tilde{x}_1)\big(\alpha_1^2\psi'^2(\tilde{x}_1) + \alpha_2^2\psi'^2(\tilde{x}_2)\big)}W(z) -\\
				&-\frac{\alpha_1\alpha_2\psi'^2(\tilde{x}_1)\psi''(\tilde{x}_1)\psi''(\tilde{x}_2) + \alpha_2^2\psi'^2(\tilde{x}_2)\big(3\psi''(\tilde{x}_1)-\psi(\tilde{x}_1)\psi'''(\tilde{x}_1)\big)}{\alpha_1^2\psi'^2(\tilde{x}_1)\big(\alpha_1^2\psi'^2(\tilde{x}_1) + \alpha_2^2\psi'^2(\tilde{x}_2)\big)}U(z)
			\end{split}
		\end{cases}
	\end{equation}
	
	\newpage
	\subsection{Results of simulation}
\begin{figure}[h]
		{\includegraphics[width=3in]{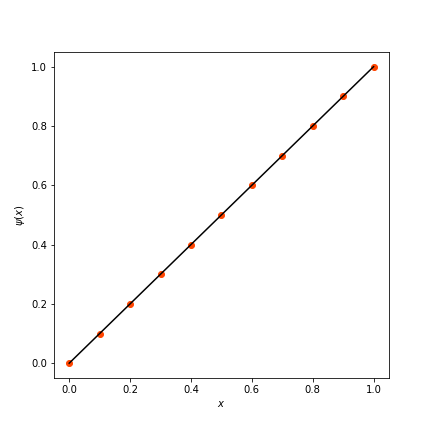}}
		\centering\caption{$\psi(x)$ for $k=1, \gamma=10$. Orange dots are points in which function $\psi(x)$ id defined. Black line is the result of linear interpolation.}
		\label{Fig3 psi_k1}
	\end{figure}

	Chosen Poisson equation is two dimensional ($n=2$), then as stated in KST $\gamma \ge 6$. We choose $\gamma = 10$ and other parameters are defined as follows: $a = \frac{1}{90}$, $\alpha_1 = 1$, $\alpha_2 = 0.10100010000000001$. For sake of simplicity, we choose $k = 1$, as for this value of $k$, inner function $\psi(x)$ is a simple identity function. Using linear interpolation, we get $\psi(x) = x$ as shown in Figure \ref{Fig3 psi_k1}.

\begin{figure}[h]
		{\includegraphics[width=4in]{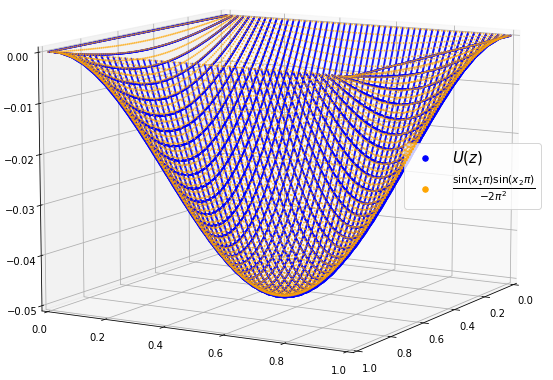}}
		\centering\caption{Solutions of Poisson equaton. Blue graph is a solution of system (\ref{System Poisson KST 0}, \ref{BC Poisson KST 0}). Orange graph is a theoretical solution of Poisson equation (\ref{PDE Poisson}, \ref{BC Poisson}).}
		\label{Fig4 PoiEq Sol}
	\end{figure}

	To solve boundary problem (\ref{System Poisson KST 0}, \ref{BC Poisson KST 0}), we use Newton-Raphson method as by introducing slight modifications it allows us to carry bifurcation analysis for boundary value problems of ordinary differential equations. The results of simulation are shown in the \textbf{Figure \ref{Fig4 PoiEq Sol}}.

	\noindent
	We can see that even using only first term in Taylor series and the simplest variant of inner function $\psi(x)$, we managed to reduce partial differential equation to the system of ordinary differential equation with the same solution.

	\section{Conclusion}
	In this paper, we introduced a new method for the transition from partial to ordinary differential equations. In comparison with existing methods, ours does not use approximation. Instead, it is based on the strict equality - Kolmogorov superposition theorem. We tested the proposed method on the Poisson equation and compared the solution of the resulting ODE system with the solution of the initial partial differential equation. In the zero approximation solutions coincided.
	
	Further, we plan to test the proposed method on non-linear partial differential equations with known solution, such as Karman equation. Secondly, we want to apply this method together with bifurcation analysis to make sure that the method does not influence the behavior of the dynamic system generated by the equation.
	\bibliographystyle{abbrv}
	\bibliography{main}

\begin{thebibliography}{10}

\bibitem{actor2018computation}
J.~Actor.
\newblock {\em Computation for the Kolmogorov Superposition Theorem}.
\newblock PhD thesis, 2018.

\bibitem{p2oReview}
J.~Awrejcewicz, V.~Krysko-Jr, L.~Kalutsky, M.~Zhigalov, and V.~Krysko.
\newblock Review of the methods of transition from partial to ordinary
  differential equations: From macro-to nano-structural dynamics.
\newblock {\em Archives of Computational Methods in Engineering}, pages 1--33,
  2021.

\bibitem{hedberg1971kolmogorov}
T.~Hedberg.
\newblock The kolmogorov superposition theorem, appendix ii to hs shapiro,
  topics in approximation theory.
\newblock {\em Lecture notes in Math}, 187:267--275, 1971.

\bibitem{khavinson1997best}
S.~I. Khavinson.
\newblock {\em Best approximation by linear superpositions (approximate
  nomography)}, volume 159.
\newblock American Mathematical Soc., 1997.

\bibitem{koppen2002training}
M.~K{\"o}ppen.
\newblock On the training of a kolmogorov network.
\newblock In {\em International Conference on Artificial Neural Networks},
  pages 474--479. Springer, 2002.

\bibitem{singhatanadgid2019kantorovich}
P.~Singhatanadgid and T.~Singhanart.
\newblock The kantorovich method applied to bending, buckling, vibration, and
  3d stress analyses of plates: A literature review.
\newblock {\em Mechanics of advanced materials and structures}, 26(2):170--188,
  2019.

\bibitem{sprecher1996numerical}
D.~Sprecher.
\newblock A numerical implementation of kolmogorov's superpositions.
\newblock {\em Neural Networks}, 9(5):765--772, 1996.

\bibitem{SPRECHER1997447}
D.~A. Sprecher.
\newblock A numerical implementation of kolmogorov's superpositions ii.
\newblock {\em Neural Networks}, 10(3):447--457, 1997.

\bibitem{sprecher2017algebra}
D.~A. Sprecher.
\newblock {\em From Algebra to Computational Algorithms: Kolmogorov and
  Hilbert's Problem 13}.
\newblock Docent Press, 2017.

\bibitem{vitushkin2004hilbert}
A.~G. Vitushkin.
\newblock On hilbert's thirteenth problem and related questions.
\newblock {\em Russian Mathematical Surveys}, 59(1):11, 2004.

\end{thebibliography}
	\nocite{SPRECHER1997447}
	\nocite{khavinson1997best}
	\nocite{sprecher2017algebra}
	\nocite{vitushkin2004hilbert}
	
	\newpage
	\section*{Appendix}
	\subsection*{A. Variational formulation of Poisson equation}\label{A}
	In order to find variational formulation of Poisson equation we start with functional (\ref{Functional Poisson Wrong}).
	\begin{equation*}\label{Functional Poisson Wrong}\tag{A.1}
		\int_{0}^{1}\int_{0}^{1}\Bigg[ \Big(\frac{\partial{u(x_1,x_2)}}{\partial{x_1}}\Big)^2 + \Big(\frac{\partial{u(x_1,x_2)}}{\partial{x_2}}\Big)^2 - 2\sin(\pi x_1)\sin(\pi x_2)u(x_1,x_2) \Bigg]dx_1dx_2 \rightarrow \underset{x_1,x_2}{extr}
	\end{equation*}
	To check whether it satisfies partial differential equation (\ref{PDE Poisson}) with boundary conditions (\ref{BC Poisson}), we solve the given problem by taking variational derivative with respect to $u(x_1, x_2)$ and equate it to zero.
	\begin{equation*}\label{Functional Poisson Wrong Variation}\tag{A.2}
		\begin{split}
			\int_{0}^{1}\int_{0}^{1}\Bigg[ 2\frac{\partial{u(x_1,x_2)}}{\partial{x_1}}\delta{u'_{x_1}(x_1,x_2)} &+ 2\frac{\partial{u(x_1,x_2)}}{\partial{x_2}}\delta{u'_{x_2}(x_1,x_2)} -\\
			&- 2\sin(\pi x_1)\sin(\pi x_2)\delta{u(x_1,x_2)} \Bigg]dx_1dx_2 = 0
		\end{split}
	\end{equation*}
	Next, we get rid of derivatives in variations using integration by parts.
	\begin{equation*}\tag{A.3}
		\begin{split}
			\int_{0}^{1}\int_{0}^{1}\Bigg[&2\nabla^2{u(x_1,x_2)}\delta{u} - 2\sin(\pi x_1)\sin(\pi x_2)\delta{u(x_1,x_2)}\Bigg]dx_1dx_2 + \\
			&+ 2\int_{0}^{1}\frac{\partial u(x_1,x_2)}{\partial x_2}\delta{u(x_1,x_2)}\Bigg|_{0}^{1}dx_1 + 2\int_{0}^{1}\frac{\partial u(x_1,x_2)}{\partial x_1}\delta{u(x_1,x_2)}\Bigg|_{0}^{1}dx_2 = 0
		\end{split}
	\end{equation*}
	Integrand satisfies initial PDE (\ref{PDE Poisson}) but natural boundary conditions do not match with (\ref{BC Poisson}). Let's add boundary terms that do not affect integrand part but change the boundary conditions so that they coincide with the PDE initial conditions (\ref{BC Poisson}). As a result, we get functional (\ref{Functional Poisson no BC}).
	\begin{equation*}\label{Functional Poisson no BC}\tag{A.4}
		\begin{split}
			\int_{0}^{1}\int_{0}^{1} \Bigg[&\Big(\frac{\partial{u(x_1,x_2)}}{\partial{x_1}}\Big)^2 + \Big(\frac{\partial{u(x_1,x_2)}}{\partial{x_2}}\Big)^2 - 2\sin(\pi x_1)\sin(\pi x_2)u(x_1, x_2)\Bigg]dx_1dx_2 -\\
			&- \int_{0}^{1}2\frac{\partial u(x_1,x_2)}{\partial x_1}u(x_1,x_2)\Bigg|_0^1dx_2 - \int_{0}^{1}2\frac{\partial u(x_1,x_2)}{\partial x_2}u(x_1,x_2)\Bigg|_0^1dx_1
		\end{split}
	\end{equation*}
	In order to apply our method for transition from partial to ordinary differential equations, we drag boundary terms under the integral. Finally, we get equivalent variational formulation for Poisson equation.
	\begin{equation*}\label{Functional Poisson}\tag{A.5}
		\begin{split}
			\int_{0}^{1}\int_{0}^{1} \bigg[&-\Big(\frac{\partial{u(x_1,x_2)}}{\partial{x_1}}\Big)^2 - \Big(\frac{\partial{u(x_1,x_2)}}{\partial{x_2}}\Big)^2 - 2\sin(\pi x_1)\sin(\pi x_2)u(x_1, x_2) -\\
			&- 2\frac{\partial^2 u(x_1,x_2)}{\partial x_1^2}u(x_1,x_2) - 2\frac{\partial^2 u(x_1,x_2)}{\partial x_2^2}u(x_1,x_2)\bigg]dx_1dx_2 \longrightarrow extr
		\end{split}
	\end{equation*}
\end{document}